\documentclass[11pt]{article}
\usepackage{amsmath}
\usepackage{amssymb}
\usepackage{listings,xcolor}
\def\a{\mathbf{a}}
\def\b{\mathbf{b}}
\def\rr{\mathbf{r}}
\def\qq{\mathbf{q}}
\def\Z{\mathbb{Z}}
\def\T{\mathcal{T}}
\def\D{\mathcal{D}}
\def\TT{\hat{\mathcal{T}}}

\newtheorem{lemma}{\hspace*{\parindent}Lemma}
\newtheorem{theorem}{\hspace*{\parindent}Theorem}
\newtheorem{proposition}{\hspace*{\parindent}Proposition}
\newtheorem{corollary}{\hspace*{\parindent}Corollary}
\title{Transformations of the hypergeometric ${}_4F_{3}$ with one unit shift: a group theoretic study}
\author{Dmitrii Karp$^{\rm a}$\footnote{Corresponding author. E-mail: D.B.\:Karp -- \emph{dmitriibkarp@tdtu.edu.vn}, E.G.\:Prilepkina --  \emph{pril-elena@yandex.ru}}~~and Elena Prilepkina$^{\rm b}$
\\[10pt]
\small{\textit{$\phantom{1}^a$Faculty of Mathematics and Statistics, Ton Duc Thang University, Ho Chi Minh City, Vietnam}}
\\\small{\textit{$\phantom{1}^b$Far Eastern Federal University and Institute of Applied Mathematics FEBRAS, Vladivostok, Russia}}}
\date{}
\begin{document}
\maketitle

\begin{abstract}
We study the group of transformations of ${}_4F_3$ hypergeometric functions evaluated at unity with one unit shift in parameters. 
We reveal the general form of this family of transformations and its group property. Next, we use explicitly known transformations  
to generate a subgroup whose structure is then thoroughly studied.  Using some known results for ${}_3F_2$ transformation 
groups, we show that this subgroup is isomorphic to the direct product of the symmetric group of degree $5$ and $5$-dimensional 
integer lattice. We investigate the relation between two-term ${}_4F_3$ transformations from our group 
and three-term ${}_3F_2$ transformations and present a method for computing the coefficients of the contiguous relations 
for ${}_3F_2$ functions evaluated at unity.  We further furnish a class of summation formulas associated with the elements of our 
group.  In the appendix to this paper, we give a collection of \textit{Wolfram Mathematica}\textsuperscript{\tiny\textregistered} 
routines facilitating the group calculations.
\end{abstract}

\bigskip


\bigskip


\bigskip

\begin{section}{Introduction and preliminaries}
Groups comprising transformation of the generalized hypergeometric functions that preserve their value at unity can 
be traced back to Kummer's formula  \cite[Corollary~3.3.5]{AAR}, see \eqref{eq:Kummer3F2} below.
These groups play an important role in mathematical physics.  In particular, the group theoretic properties of 
hypergeometric transformations constitute the key ingredient of a succinct description of the symmetries of Clebsh-Gordon's 
and Wigner's $3\!-\!j$, $6\!-\!j$ and $9\!-\!j$ coefficients from the angular momentum theory \cite{KrRao2004,Rao,RaoDobermann,RaoBook}.  
The Karlsson-Minton summation formula for the generalized hypergeometric function with integral parameter differences (IPD) 
was largely motivated by a computation of a Feymann's path integral. Furthermore, IPD hypergeometric functions appear 
in calculation of a number of integrals in high energy field theories and statistical physics \cite{ShpotSrivastava}.  
See also introduction and references in \cite{FGS} for further applications in mathematical physics and relation to Coxeter groups.

The generalized hypergeometric function \cite[2.1.2]{AAR}, \cite[Chapter~16]{NIST} is defined by the series
\begin{equation}\label{eq:pFq-defined}
{}_{p+1}F_{p}\left.\!\!\left(\begin{matrix}a_{1},\ldots,a_{p+1}\\b_{1},\ldots,b_{p}\end{matrix}\right\vert z\right)=\sum\limits_{n=0}^{\infty}\frac{(a_{1})_n\cdots(a_{p+1})_n}{n!(b_{1})_n\cdots(b_{p})_{n}}z^n
\end{equation}
whenever it converges.  When evaluated at the unit argument, $z=1$, it represents a function of $2p+1$ complex parameters with
obvious symmetry with respect to separate permutation of the $p+1$ top and the $p$ bottom parameters.
As the above series diverges at $z=1$  if the parametric excess satisfies
$\Re\left(\sum_{k=1}^{p}(b_{k}-a_{k})-a_{p+1}\right)<0$, the first problem that arises is to construct
an analytic continuation to the values of parameters in this domain.
For ${}_{3}F_{2}$ function this problem is partially solved by the transformation \cite[Corollary~3.3.5]{AAR}
\begin{equation}\label{eq:Kummer3F2}
{}_{3}F_{2}\!\left(\begin{matrix}a,b,c\\d,e\end{matrix}\right)
=\frac{\Gamma(e)\Gamma(d+e-a-b-c)}{\Gamma(e-c)\Gamma(d+e-b-a)}
{}_{3}F_{2}\!\left(\begin{matrix}d-b,d-a,c\\d,d+e-b-a\end{matrix}\right)
\end{equation}
discovered by Kummer in 1836. In the above formula we have omitted the argument $1$ from the notation
of the hypergeometric series and this convention will be adopted throughout the paper.  The series the right hand
side of \eqref{eq:Kummer3F2} converges when $\Re(e-c)>0$ so that we get the analytic continuation to this domain.
An important aspect of the above formula is that it can be applied to itself directly or after permuting
some of the top and/or bottom parameters. This leads to a family of transformations which can can be
studied by group theoretic methods.  A notable member of this family is Thomae's (1879) transformation \cite[Corollary~3.3.6]{AAR}
\begin{equation}\label{eq:Thomae3F2}
{}_{3}F_{2}\!\left(\begin{matrix}a,b,c\\d,e\end{matrix}\right)
=\frac{\Gamma(d)\Gamma(e)\Gamma(s)}{\Gamma(c)\Gamma(s+b)\Gamma(s+a)}
{}_{3}F_{2}\!\left(\begin{matrix}d-c,e-c,s\\s+a,s+b\end{matrix}\right),
\end{equation}
where $s=d+e-a-b-c$, which gave the name to the whole family of ${}_3F_{2}$ transformations generated by the
algorithm described above.  In an important work \cite{BLS1987} the authors undertook a detailed group theoretic study 
of Thomae's transformations as well as transformations for the terminating ${}_{4}F_{3}$ series and Bailey's three-term 
relations for ${}_{3}F_{2}$.  In particular, they have shown \cite[Theorem~3.2]{BLS1987} that the function
\begin{equation}\label{eq:BLS}
f(x,y,z,u,v)=\frac{{}_{3}F_{2}\!\left(\begin{matrix}x+u+v,y+u+v,z+u+v\\x+y+z+2u+v,x+y+z+u+2v\end{matrix}\right)}{\Gamma(x+y+z+2u+v)\Gamma(x+y+z+u+2v)\Gamma(x+y+z)},
\end{equation}
 is invariant with respect to the entire symmetric group $P_5$ acting on its $5$ arguments
(note that another, simpler version of this symmetry is given by  \cite[(7)]{KrRao2004}).
This symmetry was, in fact, first observed by Hardy in his 1940 lectures \cite[Notes~on~Lecture~VII]{Hardy}.
The work  \cite{BLS1987} initiated the whole stream of papers on group-theoretic  interpretations of 
hypergeometric and $q$-hypergeometric transformations. See, for instance, 
\cite{FGS,GMS,KrRao2004,Mishev,RaoDobermann,RJRJR,JeugtRao} and references therein.

We note in passing that the analytic continuation problem for general $p$ was
solved by N{\o}rlund \cite{Norlund} and Olsson \cite{Olsson} with later rediscovery by B\"{u}hring \cite{Buehring92}
without resorting to group-theoretic methods.  More recently, Kim, Rathie and Paris derived \cite[p.116]{KRP2014} 
the following transformation
\begin{equation}\label{eq:KRP2014}
{}_{4}F_3\!\!\left(\begin{matrix}a,b,c,f+1\\d,e,f\end{matrix}\right)
=\frac{\Gamma(e)\Gamma(\psi)}{\Gamma(e-c)\Gamma(\psi+c)}
{}_{4}F_3\!\!\left(\begin{matrix}d-a-1, d-b-1,c,\eta+1\\d,d+e-a-b-1,\eta\end{matrix}\right),
\end{equation}
with $\psi=d+e-a-b-c-1$ and
$$
\eta=\frac{(d-a-1)(d-b-1)f}{ab+(d-a-b-1)f}.
$$
This transformation can be iterated, but it is not immediately obvious what is the general form of the 
transformations obtained by such iterations. In our recent paper \cite[p.14, above Theorem~2]{KPGfuncmethod} we found another
identity of a similar flavor which can be viewed as a generalization of \eqref{eq:Kummer3F2}:
\begin{equation}\label{eq:KP2019}
{}_{4}F_{3}\left.\!\!\left(\!\begin{matrix}
a,b,c,f+1\\d,e,f\end{matrix}\right.\right)=
\frac{(\psi{f}-c(d-a-b))\Gamma(e)\Gamma(\psi)}{f\Gamma(e+d-a-b)\Gamma(e-c)}
{}_{4}F_{3}\left.\!\!\left(\!\begin{matrix}
d-a,d-b,c,\xi+1
\\
d,e+d-a-b,\xi\end{matrix}\right.\right),
\end{equation}
where  $\xi=f+(d-a-b)(f-c)/(e-c-1)$.  The main purpose of this paper is to present
a general form of the family of transformations of which the above two identities are particular cases,
demonstrate that this family forms a group and analyze the structure of the subgroup generated 
by explicitly known transformations \eqref{eq:KRP2014}-\eqref{eq:generator3}.
Before we delve into this analysis let us now record two more transformations
generating this subgroup.  The proof will be given in Section~4.
\begin{lemma}\label{lm:generators}
The following identities hold
\begin{equation}\label{eq:generator1}
{}_{4}F_{3}\!\!\left(\begin{matrix}a,b,c,f+1\\d,e,f\end{matrix}\right)
=\frac{(f\psi+bc)\Gamma(\psi)\Gamma(d)\Gamma(e)}{f\Gamma(a)\Gamma(\psi+b+1)\Gamma(\psi+c+1)}
{}_{4}F_{3}\!\!\left(\begin{matrix}\psi,d-a,e-a,\zeta+1\\d+e-a-c,d+e-a-b,\zeta\end{matrix}\right),
\end{equation}
where $\zeta=\psi+bc/f$, $\psi=d+e-a-b-c-1$; and
\begin{equation}\label{eq:generator3}
{}_{4}F_{3}\!\!\left(\begin{matrix}a,b,c,f+1\\d,e,f\end{matrix}\right)
=\frac{(abc+fd\psi)\Gamma(\psi)\Gamma(e)}{fd\Gamma(e-a)\Gamma(\psi+a+1)}
{}_{4}F_{3}\!\!\left(\begin{matrix}a,d-b,d-c,\nu+1\\d+1,\psi+a+1,\nu\end{matrix}\right),
\end{equation}
where $\nu=(abc+fd\psi)/(bc+f\psi)$.
\end{lemma}
Note that each ${}_4F_{3}$ function containing a parameter pair $\begin{bmatrix}f+1\\f\end{bmatrix}$ can be
decomposed into a sum of two ${}_3F_{2}$ functions (and we will demonstrate that there are numerous different decompositions of 
this type).  Hence, each of the identities \eqref{eq:KRP2014}-\eqref{eq:generator3}
can be written as a four-term relation for ${}_3F_2$.  However, it will be seen from the subsequent considerations that, in fact, 
all such relations reduce to three or even two terms, and, moreover, the structure seems to be more transparent
if we keep the ${}_4F_{3}$ function as the basic building block of our analysis. It will be revealed that the group structure 
of our transformations is closely related to that of the Thomae group generated by two-term transformations \eqref{eq:Kummer3F2} 
and \eqref{eq:Thomae3F2} and with contiguous three terms relations for ${}_3F_2$. We believe that our subgroup generated 
by \eqref{eq:KRP2014}-\eqref{eq:generator3} covers all possible two-term transformations for ${}_{4}F_{3}$ with one unit shift 
(more precisely all transformations of the form \eqref{eq:Transform} below), but we were unable to prove this 
claim and leave it as a conjecture.  

The paper is organized as follows.  In the following section we give a general form of the transformations exemplified above and 
prove that they form a group. We further demonstrate that this group is isomorphic to a subgroup of $\mathrm{SL}(\Z)$ 
(integer matrices with unit determinant).  In Section~3, we give a comprehensive analysis of  the structure of 
the subgroup generated by the transformations \eqref{eq:KRP2014}-\eqref{eq:generator3} by showing that it is 
isomorphic to a direct product of the symmetric group $P_5$ and the integer lattice $\Z^5$.   
In Section~4 we explore the relation between our transformations and three-term relations for ${}_3F_2$ hypergeometric function. 
In particular, we show that the contiguous relations for ${}_3F_2$ functions studied recently in \cite{EbisuIwasaki} 
can also be computed from the elements of our group.  Section~5 contains the proof of Lemma~\ref{lm:generators} 
and the Appendix contains explicit forms of some key elements of our subgroup and several  
\textit{Wolfram Mathematica}\textsuperscript{\tiny\textregistered}  routines facilitating the group calculations.
\end{section}

\begin{section}{The group structure of the unit shift ${}_4F_{3}$ transformations}
Inspecting the ${}_4F_{3}$ transformations presented in Section~1 we see that they share a common structure 
that we will present below. To this end, let $\rr=(a,b,c,d,e,1)^T$ be the column vector and define
\begin{equation}\label{eq:4F3unitshift}
F(\rr,f)={}_{4}F_{3}\!\left(\begin{matrix}a,b,c,f+1\\d,e,f\end{matrix}\right).
\end{equation}
All transformations found in Section~1 have the following general form
\begin{equation}\label{eq:Transform}
F(\rr,f)=C(\rr,f)F(D\rr,\eta),
\end{equation}
where $D$ is a unit determinant $6\times 6$ matrix with integer entries and the bottom row $(0,0,0,0,0,1)$;
\begin{equation}\label{eq:Transeta}
\eta=\frac{\varepsilon{f}+\lambda(\rr)}{\alpha(\rr)f+\beta(\rr)},
\end{equation}
where $\varepsilon\in\{0,1\}1$, $\lambda(\rr)$, $\alpha(\rr)$ and $\beta(\rr)$ are rational functions of the arguments $a,b,c,d,e$ 
(some of them may vanish identically, but $\lambda=1$ if $\varepsilon=0$). The coefficient $C(\rr,f)$ has the form
\begin{equation}\label{eq:Cstructure}
C(\rr,f)=\frac{N(\rr)f+P(\rr)}{K(\rr)f+L(\rr)}
\end{equation}
where $N(\rr)$, $P(\rr)$, $K(\rr)$, $L(\rr)$ are (possibly vanish) functions of $\Gamma$-type 
by which we mean ratios of products of gamma functions whose arguments are integer linear combinations of the 
components of $(a,b,c,d,e,1)$.  When $N(\rr)\ne0$  we will additionally require that the 
ratio  $P(\rr)/N(\rr)$ be a rational function of parameters. In fact, this last requirements is redundant, 
but in order to  avoid it the following claim is needed: 
the ratio $F_2(\rr)/F_1(\rr)$ with $F_i,$ $i=1,2$, defined in \eqref{eq:part1}, 
is not a function of gamma type for general parameters. We were unable to find a proof of this claim in the literature 
although it seems to be generally accepted to be true.

Formula \eqref{eq:Transform} defines a transformation $T$ characterized by the matrix $D$  and the functions 
$C(\rr,f)$, $\eta=\eta(\rr,f)$. Two  such transformations $T_1$, $T_2$ will be considered equal 
if $D_1=D_2$, $C_1(\rr,f)\equiv C_2(\rr,f)$ and $\eta_1(\rr,f)\equiv \eta_2(\rr,f)$.

According to the elementary relation $(f+1)_n=(f)_n(1+n/f)$, we have
\begin{equation}\label{eq:partition}
F(\rr,f)=F_1(\rr)+\frac{1}{f}F_2(\rr),
\end{equation}
where
\begin{equation}\label{eq:part1}
F_1(\rr)={}_{3}F_{2}\!\left(\begin{matrix}a,b,c\\d,e\end{matrix}\right),~~~  F_2(\rr)=\frac{abc}{de}{}_{3}F_{2}\!\left(\begin{matrix}a+1,b+1,c+1\\d+1,e+1\end{matrix}\right).
\end{equation}

It is not immediately obvious if the composition of two transformations \eqref{eq:Transform} with $\eta$ and $C$ having
the forms \eqref{eq:Transeta} and \eqref{eq:Cstructure}, respectively, should have the same form.
The following theorem shows that it is indeed the case and these transformations form a group.
\begin{theorem}\label{th:Group}
Each transformation  \eqref{eq:Transform} necessarily has the form
\begin{equation}\label{eq:generalform}
F(\rr,f)=M(\rr)\frac{\varepsilon{f}+\lambda(\rr)}{f}F(D\rr,\eta),~~\text{where}~~\eta=\frac{\varepsilon{f}+\lambda(\rr)}{\alpha(\rr)f+\beta(\rr)},
\end{equation}
$M(\rr)$ is a function of $\Gamma$-type, $\varepsilon\in\{0,1\}$, $\lambda(\rr)$, $\alpha(\rr)$, $\beta(\rr)$ are rational functions of the arguments $a,b,c,d,e$
\emph{(}possibly vanishing but with $\lambda=1$ if $\varepsilon=0$\emph{)}.

\smallskip

The collection $\T$ of transformations \eqref{eq:generalform} forms a group with respect to composition.
More explicitly, if $T_1,T_2\in\T$ with parameters indexed correspondingly, then $T=T_2\circ{T_1}$ is given by

\bigskip

\emph{(I)} If $\varepsilon_{1}\varepsilon_{2}+\alpha_1(\rr)\lambda_2(D_1\rr)\ne0$, then $\varepsilon=1$,
$M(\rr)=M_1(\rr)M_2(D_1\rr)(\varepsilon_{1}\varepsilon_{2}+\alpha_1(\rr)\lambda_2(D_1\rr))$,
$$
\lambda(\rr)=\frac{\varepsilon_{2}\lambda_1(\rr)+\lambda_2(D_{1}\rr)\beta_{1}(\rr)}
{\varepsilon_{1}\varepsilon_{2}+\alpha_1(\rr)\lambda_2(D_1\rr)},
~~~\alpha(\rr)=\frac{\varepsilon_1\alpha_2(D_{1}\rr)+\alpha_{1}(\rr)\beta_{2}(D_{1}\rr)}
{\varepsilon_{1}\varepsilon_{2}+\alpha_1(\rr)\lambda_2(D_1\rr)},
$$
$$
\beta(\rr)=\frac{\lambda_{1}(\rr)\alpha_{2}(D_{1}\rr)+\beta_{1}(\rr)\beta_{2}(D_1\rr)}
{\varepsilon_{1}\varepsilon_{2}+\alpha_1(\rr)\lambda_2(D_1\rr)},~~~~D=D_{2}D_{1}.
$$

\medskip

\emph{(II)} If $\varepsilon_{1}\varepsilon_{2}+\alpha_1(\rr)\lambda_2(D_1\rr)=0$, then  $\varepsilon=0$,
$M(\rr)=M_1(\rr)M_2(D_1\rr)(\varepsilon_{2}\lambda_1(\rr)+\lambda_2(D_{1}\rr)\beta_{1}(\rr))$,
$$
\lambda(\rr)=1,
~~~\alpha(\rr)=\frac{\varepsilon_1\alpha_2(D_{1}\rr)+\alpha_{1}(\rr)\beta_{2}(D_{1}\rr)}{\varepsilon_{2}\lambda_1(\rr)+\lambda_2(D_{1}\rr)\beta_{1}(\rr)},
$$
$$
\beta(\rr)=\frac{\lambda_{1}(\rr)\alpha_{2}(D_{1}\rr)+\beta_{1}(\rr)\beta_{2}(D_1\rr)}{\varepsilon_{2}\lambda_1(\rr)+\lambda_2(D_{1}\rr)\beta_{1}(\rr)},
~~~D=D_{2}D_{1}.
$$

Each $T\in\T$ of the form \eqref{eq:generalform} has an inverse $T^{-1}$ determined  by the 
parameters $\hat{\varepsilon}$, $\hat{M}(\rr)$,  $\hat{\lambda}(\rr),$ $\hat{\alpha}(\rr)$, $\hat{\beta}(\rr)$, $\hat{D}$ 
given by:

\emph{(III)} If $\beta(\rr)\ne0$, then $\hat{\varepsilon}=1$ and
$$
\hat{M}(\rr)=\frac{\beta(D^{-1}\rr)}{M(D^{-1}\rr)(\varepsilon\beta(D^{-1}\rr)-\alpha(D^{-1}\rr)\lambda(D^{-1}\rr)},
~~~\hat{\lambda}(\rr)=-\frac{\lambda(D^{-1}\rr)}{\beta(D^{-1}\rr)},
$$
$$
\hat{\alpha}(\rr)=-\frac{\alpha(D^{-1}\rr)}{\beta(D^{-1}\rr)},~~~
\hat{\beta}(\rr)=\frac{\varepsilon}{\beta(D^{-1}\rr)},~~~\hat{D}=D^{-1}.
$$

\emph{(IV)} If $\beta(\rr)=0$, then $\hat{\varepsilon}=0$ and
$$
\hat{M}(\rr)=\frac{1}{M(D^{-1}\rr)\alpha(D^{-1}\rr)},
~~
\hat{\lambda}(\rr)=1,~~ \hat{\alpha}(\rr)=\frac{\alpha(D^{-1}\rr)}{\lambda(D^{-1}\rr)},~~
\hat{\beta}(\rr)=-\frac{\varepsilon}{\lambda(D^{-1}\rr)},~~~\hat{D}=D^{-1}.
$$
\end{theorem}

\textbf{Proof.} We start by showing that the form of the coefficient
$C(\rr,f)=(Nf+P)/(Kf+L)$ defined in \eqref{eq:Cstructure} is restricted to 
\begin{equation}\label{eq:koef}
C(\rr,f)=M+W/f,
\end{equation}
where $M=M(\rr)$, $W=W(\rr)$ are some functions of $\Gamma$-type, possibly one of them vanishing. It follows 
from \eqref{eq:Cstructure} and \eqref{eq:partition} that  transformation \eqref{eq:Transform} is equivalent to
\begin{equation}\label{eq:eq1}
\frac{F_1(\rr)f+F_2(\rr)}{f}=\frac{(Nf+P)(F_1(D\rr)\eta+F_2(D\rr))}{(Kf+L)\eta},
\end{equation}
where $N=N(\rr),$ $P=P(\rr),$ $K=K(\rr),$ $L=L(\rr)$. Solving this equation we get
$$
\eta=\frac{f(fN+P)F_2(D\rr)}{LF_2(\rr)+fKF_2(\rr)-f^2NF_1(D\rr)-fPF_1(D\rr)+f^2KF_1(\rr)+fLF_1(\rr)}.
$$
In order that  $\eta$ had the form \eqref{eq:Transeta} the following identity must hold
\begin{multline}\label{eq:yaqq}
f(fN+P)F_2(D\rr)(\alpha{f}+\beta)
\\
=(\varepsilon{f}+\lambda)({LF_2(\rr)+fKF_2(\rr)-f^2NF_1(D\rr)-fPF_1(D\rr)+f^2KF_1(\rr)+fLF_1(\rr)}).
\end{multline}
The free term of the cubic on the right hand side equals $\lambda{L}F_2(\rr)$ while it vanishes on the left hand side, 
so that  $\lambda{L}=0$. If  $L=0$ we obtain  \eqref{eq:koef}. Otherwise, if $\lambda=0$ identity \eqref{eq:yaqq} takes the form
\begin{multline}\label{eq:yaqq1}
(fN+P)F_2(D\rr)(\alpha{f}+\beta)
\\
=LF_2(\rr)+fKF_2(\rr)-f^2NF_1(D\rr)-fPF_1(D\rr)+f^2KF_1(\rr)+fLF_1(\rr).
\end{multline}
If $N=0$,  then $K=0$ and we again arrive at \eqref{eq:koef}. If $N\ne0$ the value $f=-P/N$ must be a root of the quadratic 
on the right hand side of \eqref{eq:yaqq1}.  In other words, we must have
$$
LF_2(\rr)-\frac{P}{N}KF_2(\rr)-\frac{P^2}{N^2}NF_1(D\rr)+\frac{P}{N}PF_1(D\rr)+\frac{P^2}{N^2}KF_1(\rr)-\frac{P}{N}LF_1(\rr)=0
$$
or
$$
\left(L-\frac{P}{N}K\right)\left(F_2(\rr)-\frac{P}{N}F_1(\rr)\right)=0.
$$
Equality $L=PK/N$ again leads to \eqref{eq:koef}. The equality $F_2(\rr)= PF_1(\rr)/N$ is impossible for rational ${P}/{N}$, 
as demonstrated by Ebisu and Iwasaki in \cite[Theorem~1.1]{EbisuIwasaki} which proves our claim \eqref{eq:koef}. 
If ${P}/{N}$ is a function of gamma type then so is $F_2(\rr)/F_1(\rr)$ which would contradict the claim made 
before the theorem, but as we could not find a proof of this claim we explicitly prohibit this situation 
in the definition of $C(\rr,f)$.

Substituting  $(Nf+P)/(Kf+L)$ by  $M+W/f$ in \eqref{eq:eq1} we can now express $\eta$ as follows:
\begin{equation}\label{eq:relation}
\eta=-\frac{(Mf+W)F_2(D\rr)}{(MF_1(D\rr)-F_1(\rr))f+F_1(D\rr)W-F_2(\rr)}.
\end{equation}

Next suppose $M\ne0$. Then $C(\rr,f)=M(\varepsilon{f}+W/M)/f$ with $\varepsilon=1$.  Comparison of \eqref{eq:relation} 
with \eqref{eq:Transeta} yields $W/M=\lambda$  which proves that the transformation \eqref{eq:Transform} must have 
the form \eqref{eq:generalform}.  Moreover,
$$
\alpha=-\frac{M\varepsilon F_1(D\rr)-F_1(\rr)}{MF_2(D\rr)},
~~\beta=-\frac{F_1(D\rr)\lambda M-F_2(\rr)}{MF_2(D\rr)}.
$$
These equalities can be rewritten as the system
 \begin{equation}\label{eq:sistem}
\left\{\begin{array}{l}
F_1(\rr)=M(\varepsilon F_1(D\rr)+\alpha {F_2({D\rr})}\bf),\\
F_2(\rr)=M(\lambda F_1(D\rr) +\beta F_2({D\rr})).\end{array}\right.
\end{equation}

Suppose now that $M=0$,  $W\ne0$. Then $C(\rr,f)=W(\varepsilon f+\lambda)/f$
with $\varepsilon=0$,  $\lambda=1$.  From  \eqref{eq:relation} we have
$$
\eta=\frac{\varepsilon f+\lambda}{\alpha f+\beta},
$$
where again $\varepsilon=0$,  $\lambda=1$, and $\alpha=F_1(\rr)/(WF_2(D\rr)),~~~\beta=-(WF_1(D\rr)-F_2(\rr))/(WF_2(D\rr))$ or
\begin{equation}\label{eq:sistem2}
\left\{\begin{array}{l}
F_1(\rr)= W{\alpha}F_2(D\rr)=W (\varepsilon F_1(D\rr)+\alpha {F_2(D\rr)}\bf),
\\
F_2(\rr)=W(F_1(D\rr) +\beta F_2(D\rr))=W(\lambda F_1(D\rr) +\beta F_2(D\rr)).
\end{array}\right.
\end{equation}
Renaming $W$ into $M$ we have thus proved that the transformation again has the form \eqref{eq:generalform} 
and the system  \eqref{eq:sistem} is satisfied.

The computation of composition is straightforward:
\begin{multline*}
T=T_2\circ{T_1}  \Longleftrightarrow  T: F(\rr,f)=M_1(\rr)\frac{\varepsilon_{1}f+\lambda_1(\rr)}{f}F(D_1\rr,\eta_1)
\\
=M_1(\rr)M_2(D_1\rr)\frac{\varepsilon_{1}f+\lambda_1(\rr)}{f}
\frac{\varepsilon_{2}(\varepsilon_{1}f+\lambda_1(\rr))/(\alpha_1(\rr)f+\beta_1(\rr))+\lambda_2(D_1\rr)}
{(\varepsilon_{1}f+\lambda_1(\rr))/(\alpha_1(\rr)f+\beta_1(\rr))}F(D_2D_1\rr,\eta_2)
\\
=M_1(\rr)M_2(D_1\rr)
\frac{\varepsilon_{2}(\varepsilon_{1}f+\lambda_1(\rr))+\lambda_2(D_1\rr)({\alpha_1(\rr)f+\beta_1(\rr)})}{f}
F(D_2D_1\rr,\eta_2)
\\
=M_1(\rr)M_2(D_1\rr)
\frac{[\varepsilon_1\varepsilon_2+\alpha_1(\rr)\lambda_2(D_1\rr)]f+\varepsilon_{2}\lambda_1(\rr)+\beta_{1}(\rr)\lambda_2(D_1\rr)}{f}
F(D_2D_1\rr,\eta).
\end{multline*}
If $\varepsilon_{1}\varepsilon_{2}+\alpha_1(\rr)\lambda_2(D_1\rr)\ne0$, we can divide by this quantity leading to case (I).
If it vanishes we get case (II).
Given a transformation $T\in\T$ of the from \eqref{eq:generalform} it is rather straightforward to compute its inverse.
We omit the details.
$\hfill\square$.

\textbf{Remark}. Theorem~\ref{th:Group} implies that each transformation $t\in\T$ is uniquely characterized by the collection
$\{\varepsilon, M(\rr), \lambda(\rr), \alpha(\rr), \beta(\rr), D\}$, where $\varepsilon\in\{0,1\}$, $M(\rr)$ is a function of
gamma type, $\lambda(\rr)$, $\alpha(\rr)$ and $\beta(\rr)$ are rational functions of parameters $a,b,c,d,e$ and $D$ is $6\times6$
unit determinant integer matrix with bottom row $(0,\ldots,0,1)$.  We will express this fact by writing
$T\sim\{\varepsilon, M(\rr), \lambda(\rr), \alpha(\rr), \beta(\rr), D\}$.  Occasionally, we will omit the dependence on $\rr$ in 
the notation of the functions $M(\rr)$, $\lambda(\rr)$, $\alpha(\rr)$, $\beta(\rr)$ for brevity.

Note that  for $\varepsilon=1$ and non-vanishing  $\alpha$, $\beta$ and $\lambda$ the system \eqref{eq:sistem} 
takes the form of ${}_4F_{3}\to{}_3F_{2}$ reduction formulas
\begin{equation*}
\left\{\begin{array}{l}
F(D\rr,\alpha(\rr)^{-1})=M(\rr)^{-1}F_1(\rr),
\\
F(D\rr,\lambda(\rr)/\beta(\rr))=(M(\rr)\lambda(\rr))^{-1}F_2(\rr).
\end{array}\right.
\end{equation*}

Next, we clarify the structure of the group $\T$ further.  The composition rule  involves all parameters
$M(\rr)$, $\lambda(\rr)$, $\alpha(\rr)$, $\beta(\rr)$ and $D$. The following theorem implies that the matrix $D$
determines all other parameters uniquely.  Denote by $\widehat{\mathrm{SL}}(n,\Z)$ the subgroup of the special 
linear group $\mathrm{SL}(n,\Z)$ of $n\times{n}$ integer matrices with unit determinant comprising matrices whose 
last row  has the form $(0,\ldots,0,1)$.

\begin{theorem}\label{th:matrix-represent}
The mapping $T\sim\{\varepsilon, M(\rr), \lambda(\rr), \alpha(\rr), \beta(\rr), D_T\}\to{D_T}$ is isomorphism, so
that the group $(\T,\circ)$ is isomorphic to a subgroup  of $\widehat{\mathrm{SL}}(n,\Z)$ which we denote by $(\mathcal{D}_{\T}, \cdot)$.
\end{theorem}

\textbf{Proof.} One direction is clear: each transformation $T\in\T$ by construction defines a matrix 
$D_T\in\widehat{\mathrm{SL}}(n,\Z)$ and the composition rule (I), (II) in Theorem~\ref{th:Group} involves the 
product of matrices. Hence, to establish our claim it remains to prove that the kernel of the homomorphism $\T\to D_T$ is trivial.
Assume the opposite: there exists a transformation $T\in\T$ with the identity matrix $D\rr=\rr$ and non-trivial parameters
$\varepsilon$, $M$, $\lambda$, $\alpha$, $\beta$.  The system \eqref{eq:sistem} then takes the form
\begin{equation}\label{eq:sistemId}
\left\{\begin{array}{l}
(1-M\varepsilon)F_1(\rr)= M\alpha {F_2({\rr})},
\\
M\lambda F_1(\rr)=(1-M\beta)F_2({\rr}).
\end{array}\right.
\end{equation}
If $\alpha=\lambda=0$ we get $M=\varepsilon=1$ from the first equation and $\beta=1$ from the second equation
which amounts to the trivial identity transformation.  We will show that all other cases are impossible.
Indeed, Ebisu and Iwasaki demonstrated in \cite[Theorem~1.1]{EbisuIwasaki} that the functions
$F_1(\rr)$ and  $F_2(\rr)$  are linearly independent over the field of rational functions of parameters.
If $\alpha=0$  and $\lambda\ne0$,  then $M=\varepsilon=1$ from the first equation and
$F_1(\rr)/F_2(\rr)=(1-\beta)/\lambda$ from the second equation contradicting linear independence.
Similarly, if $\alpha\ne0$ and $\lambda=0$, then $M=1/\beta$ from the second equation, so that
$F_2(\rr)/F_1(\rr)=(1-\varepsilon/\beta)/(\alpha/\beta)$ is rational from the first equation leading again to contradiction.
Finally, if both $\alpha\ne0$ and $\lambda\ne0$ we arrive at the identities
$$
\frac{F_2(\rr)}{F_1(\rr)}=\frac{1-M\varepsilon}{\alpha M}=\frac{\lambda M}{1-M\beta}~~\Rightarrow~~(1-M\beta)(1-M\varepsilon)=\alpha{\lambda}M^2.
$$
Linear independence of the functions  $F_1(\rr)$, $F_2(\rr)$ over rational functions implies that the function $M=M(\rr)$ must 
be a ratio of products of gamma functions  irreducible to a rational function. On the other hand, by the ultimate 
equality $M(\rr)$ solves the quadratic equation with rational coefficients:
$$
M=M(\rr)=\mu(\rr)\pm\sqrt{\nu(\rr)}
$$
with rational $\mu$, $\nu$. It is easy to see that this is not possible as $\Gamma$ is meromorphic with infinite number of poles
and no branch points, while  $\mu(\rr)\pm\sqrt{\nu(\rr)}$ may only have a finite number of poles and zeros and has branch points.
$\hfill\square$

\end{section}

\begin{section}{The subgroup of $\T$ generated by known transformations}

We can now rewrite the transformations \eqref{eq:KRP2014}-\eqref{eq:generator3} in the standard form \eqref{eq:generalform}. 
Denote by  $\psi=d+e-a-b-c-1$ the parametric excess of the function on the left hand side of \eqref{eq:generalform}.  
Identity \eqref{eq:generator1} is determined by the following set of parameters
\begin{subequations}\label{eq:T1}
\begin{align}
&M_1=\frac{\Gamma(\psi+1)\Gamma(d)\Gamma(e)}{\Gamma(a)\Gamma(d+e-a-c)\Gamma(d+e-a-b)},~~~
\varepsilon_1=1,~~~\lambda_1=\frac{bc}{\psi},
\\[5pt]
&D_1=\begin{bmatrix}
-1&-1&-1&1&1&-1\\
-1&0&0&1&0&0\\
-1&0&0&0&1&0\\
-1&0&-1&1&1&0\\
-1&-1&0&1&1&0\\
0&0&0&0&0&1
\end{bmatrix},~~~\alpha_1=\frac{1}{\psi},~~~ \beta_1=0.
\end{align}
\end{subequations}
We will call this transformation $T_1$.

The standard form \eqref{eq:generalform} of identity \eqref{eq:KP2019} is characterized by the following parameters:
\begin{subequations}\label{eq:T2}
\begin{align}
&M_2=\frac{\Gamma(e)\Gamma(\psi+1)}{\Gamma(e+d-a-b)\Gamma(e-c)},~~~
\varepsilon_2=1,~~~\lambda_2=\frac{c(-d+a+b)}{\psi},
\\[5pt]\label{eq:T2-2}
&D_2=\begin{bmatrix}
-1&	0&	0&	1&	0&	0\\
0&	-1&	0&	1&	0&	0\\
0&	0&	1&	0&	0&	0\\
0&	0&	0&	1&	0&	0\\
-1&	-1&	0&	1&	1&	0\\
0&	0&	0&	0&	0&	1\\
\end{bmatrix},~~~\alpha_2=0,~~~\beta_2=\frac{e-c-1}{\psi}.
\end{align}
\end{subequations}
We will call this transformation $T_2$.

The standard parameters of transformation \eqref{eq:generator3} are given by
\begin{subequations}\label{eq:T3}
\begin{align}
&M_3=\frac{\Gamma(\psi+1)\Gamma(e)}{\Gamma(e-a)\Gamma(e+d-b-c)},~~~\varepsilon_3=1,
~~~\lambda_3=\frac{abc}{d\psi},
\\[5pt]
&D_3=\begin{bmatrix}
1&	0&	0&	0&	0&	0\\
0&	-1&	0&	1&	0&	0\\
0&	0&	-1&	1&	0&	0\\
0&	0&	0&	1&	0&	1\\
0&	-1&	-1&	1&	1&	0\\
0&	0&	0&	0&	0&	1
\end{bmatrix},~~~~\alpha_3=\frac{1}{d},~~~\beta_3=\frac{bc}{d\psi}.
\end{align}
\end{subequations}
We will call this transformation $T_3$.

Finally, transformation \eqref{eq:KRP2014} in the standard form \eqref{eq:generalform} is parameterized by
\begin{subequations}\label{eq:T4}
\begin{align}
&M_4=\frac{\Gamma(e)\Gamma(\psi)}{\Gamma(e-c)\Gamma(\psi+c)},~~~
\varepsilon_4=1,~~~\lambda_4=0,~~~\alpha_4=\frac{d-a-b-1}{(d-a-1)(d-b-1)},
\\[5pt]
&D_4=\begin{bmatrix}
-1&	0&	0&	1&	0&	-1\\
0&	-1&	0&	1&	0&	-1\\
0&	0&	1&	0&	0&	0\\
0&	0&	0&	1&	0&	0\\
-1&	-1&	0&	1&	1&	-1\\
0&	0&	0&	0&	0&	1\\
\end{bmatrix},~~~\beta_4=\frac{ab}{(d-a-1)(d-b-1)}.
\end{align}
\end{subequations}
We will call this transformation $T_4$. It is easy to see that it is of order $2$, i.e. $T_4^2=I$.

The four transformations $T_1$, $T_2$, $T_3$, $T_4$ (or, equivalently, \eqref{eq:KRP2014}, \eqref{eq:KP2019},
\eqref{eq:generator1} and \eqref{eq:generator3}) combined with permutations of the upper and lower parameters 
generate a subgroup of $\T$ which we will call $\TT$.  Isomorphism established in Theorem~\ref{th:matrix-represent} 
induces an isomorphism between $\TT$ and a subgroup of $\widehat{\mathrm{SL}}(n,\Z)$ which we denote by $\D_{\TT}$.

A complete characterization of $\TT$ and $\D_{\TT}$ will follow. Before we turn to it,
we remark that to our belief, the complete group $\T$ contains no elements other than those in $\TT$.
We were unable, however, to prove this claim.  Let us thus state it as a conjecture.

\textbf{Conjecture.}  The subgroup $\TT$ generated by the transformations \eqref{eq:T1}, \eqref{eq:T2}, \eqref{eq:T3} and
\eqref{eq:T4} coincides with the entire group $\T$ of all transformations of the form \eqref{eq:Transform} or, equivalently,
of the form \eqref{eq:generalform}.

Denote by $S_{j}$, $j=1,\ldots,5$, the transformation shifting the $j$-th component of the parameter vector $\rr$ by $+1$,
i.e. $S_{j}$ is characterized by the matrix $\hat{S}_{j}$ such that $\hat{S}_{1}\rr=(a+1,b,c,d,e,1)$,
$\hat{S}_{2}\rr=(a,b+1,c,d,e,1)$, etc.  It is not \emph{a priori} obvious that such transformations should 
exist among the elements of $\TT$. The following theorem shows that it is indeed the case.

\begin{theorem}\label{th:shift}
The group $\TT$ contains the transformations $S_{j}$, $j=1,\ldots,5$.
\end{theorem}
\textbf{Proof}. Due to permutation symmetry it is clearly sufficient to display the transformations $S_1$ and $S_4$.
We will need the inverse of the transformation $T_1$ defined in \eqref{eq:T1}. Using Theorem~\ref{th:Group} we calculate
\begin{equation}\label{eq:T1inverse}
{}_{4}F_3\!\!\left(\begin{matrix}a,b,c,f+1\\d,e,f\end{matrix}\right)
=\frac{\hat{M}_1}{f}
{}_{4}F_3\!\!\left(\begin{matrix}d+e-a-b-c-1,d-a-1,e-a-1,\hat{\eta}_1+1\\d+e-a-c-1,d+e-a-b-1,\hat{\eta}_1\end{matrix}\right),
\end{equation}
where
$$
\hat{M}_1=\frac{\Gamma(d)\Gamma(e)\Gamma(\psi)}{\Gamma(\psi+b)\Gamma(\psi+c)\Gamma(a)},~~~\hat{\varepsilon}_1=0,~~~\hat{\lambda}_1=1,~~\hat{\alpha}_1=\frac{1}{(d-a-1)(e-a-1)},
$$
$$
\hat{\beta}_1=\frac{-a}{(d-a-1)(e-a-1)},~~\text{so that}~~\hat{\eta}_1=\frac{(d-a-1)(e-a-1)}{f-a}.
$$
Next, exchanging the roles of $d+e-a-b-1$ and $d$ and the roles of $d-a-1$ and $c$ in \eqref{eq:KRP2014} or, equivalently, 
post-composing $T_4$ with permutation $(13)(45)$ we will obtain a transformation that we call $\hat{T}_4$. 
Then $\hat{T}_4\circ{\hat{T}_4}$ takes the form
\begin{equation}\label{eq:iteratedMRP}
{}_{4}F_3\!\!\left(\begin{matrix}a,b,c,f+1\\d,e,f\end{matrix}\right)
=\frac{\Gamma(e)\Gamma(d)\Gamma(\psi)}{\Gamma(b+\psi)\Gamma(c+\psi)\Gamma(a+1)}
{}_{4}F_3\!\!\left(\begin{matrix}\psi-1,e-a-1,d-a-1,\tilde{\eta}_4+1\\c+\psi,b+\psi,\tilde{\eta}_4\end{matrix}\right)
\end{equation}
with
$$
\tilde{\eta}_4=\frac{(\psi-1)(e-a-1)(d-a-1)f}{abc+(1+2a+a^2-bc-d-ad-e-ae+de)f}.
$$
Applying $T_1^{-1}$ to the right hand side of \eqref{eq:iteratedMRP} we obtain the transformation $S_1$:
\begin{equation}\label{eq:killuppershift}
{}_{4}F_3\!\!\left(\begin{matrix}a,b,c,f+1\\d,e,f\end{matrix}\right)
=M\frac{\varepsilon{f}+\lambda}{f}
{}_{4}F_3\!\!\left(\begin{matrix}a+1,b,c,\eta+1\\d,e,\eta\end{matrix}\right),
\end{equation}
where $\varepsilon=1$, and
$$
M=1-\frac{bc}{(d-a-1)(e-a-1)},~~~\lambda=\frac{abc}{a^2-bc+(d-1)(e-1)-a(d+e-2)},
$$
$$
\alpha=\frac{d+e-a-b-c-2}{a^2-bc+(d-1)(e-1)-a(d+e-2)},~~\beta=-\frac{a(d+e-a-b-c-2)}{a^2-bc+(d-1)(e-1)-a(d+e-2)}.
$$
According to \eqref{eq:generalform} we thus obtain the following expression for $\eta$
$$
\eta=\frac{abc+(1+2a+a^2-bc-d-ad-e-ae+de)f}{a(2+a+b+c-d-e)-(2+a+b+c-d-e)f}.
$$

Application of the transformation $T_3$ given by \eqref{eq:T3} to itself yields $T_3\circ{T_3}$ in the form:
$$
{}_{4}F_{3}\!\!\left(\begin{matrix}a,b,c,f+1\\d,e,f\end{matrix}\right)
=\frac{a(d-b)(d-c)(bc+f\psi)+(d+1)(e-a)(abc+fd\psi)}{fd(d+1)e\psi}
{}_{4}F_{3}\!\!\left(\begin{matrix}a,b+1,c+1,\tilde{\eta}_3+1\\d+2,e+1,\tilde{\eta}_3\end{matrix}\right),
$$
where
$$
\tilde{\eta}_3=\frac{a(d-b)(d-c)(bc+f\psi)+(d+1)(e-a)(abc+fd\psi)}{(d-b)(d-c)(bc+f\psi)+(e-a)(abc+fd\psi)}.
$$
On the other hand, using \eqref{eq:T1inverse} we compute $T_{1}^{-2}$ as follows:
$$
{}_{4}F_3\!\!\left(\begin{matrix}a,b,c,f+1\\d,e,f\end{matrix}\right)
=\frac{(d-1)(e-1)(f-a)}{f(d-a-1)(e-a-1)}
{}_{4}F_3\!\!\left(\begin{matrix}a,b-1,c-1,\hat{\eta}_1'+1\\d-1,e-1,\hat{\eta}_1'\end{matrix}\right)
$$
with
$$
\hat{\eta}_1'=\frac{(b-1)(c-1)(f-a)}{(d-a-1)(e-a-1)-\psi(f-a)}.
$$
Comparing these formulas we see that the composition $T_1^{-2}\circ{T_3^2}$ gives the transformation $S_4$ shifting $d\to{d+1}$ while $a,b,c,e$ remain intact:
\begin{equation}\label{eq:killbottom}
{}_{4}F_3\!\!\left(\begin{matrix}a,b,c,f+1\\d,e,f\end{matrix}\right)
=\frac{f+\lambda}{f}
{}_{4}F_3\!\!\left(\begin{matrix}a,b,c,\eta+1\\d+1,e,\eta\end{matrix}\right),
\end{equation}
so that $\varepsilon=1$, $M=1$,
$$
\lambda=\frac{abc}{d(d+e-a-b-c-1)},~~\alpha=\frac{1}{d},
~~\beta=\frac{(b-d)(c-d)+a(b+c-d)}{d(d+e-a-b-c-1)},~~\eta=\frac{\varepsilon{f}+\lambda}{\alpha{f}+\beta}.~~~~\square
$$

Each transformation $S_j$, $j=1,\ldots,5$, obviously generates a subgroup of $\TT$ isomorphic to
$\Z$ - the additive group of integers. Hence, in the parlance of group theory, the above theorem
can be restated and enhanced as follows.

\begin{corollary}
The group $\TT$ contains a subgroup $\mathcal{S}$ isomorphic to the $5$-dimensional integer lattice $\Z^5$.
Furthermore, this subgroup is normal.
\end{corollary}
\textbf{Proof.}  By the previous theorem we only need to prove normality. Denote by $\mathcal{S}$ the subgroup of the matrix group $\D_{\TT}$ generated by the shift matrices $\hat{S}_{j}$, $j=1,\ldots,5$.  Clearly, $\mathcal{S}$ comprises $6\times6$ matrices whose principal $5\times5$ sub-matrix equals the identity matrix $I_5$,  the $6$-th row is $(0,\dots,0,1)$ and the $6$-th column is  $(k_1,\ldots,k_5,1)$ for some $k_i\in\Z$.
As all elements of $\D_{\TT}$ have integer entries and the bottom row $(0,\dots,0,1)$
it is easy to see that for any shift matrix $S\in\mathcal{S}$ and any matrix $D\in\D_{\TT}$ both products
$DS$ and $SD$ have the principal $5\times5$ sub-matrix equal to that of $D$ and the last column of the form
$(k_1,\ldots,k_5,1)$ for some $k_i\in\Z$. Running over all elements of $\mathcal{S}$ while keeping $D$ fixed we
see that the left and right conjugacy classes of the element $D$ with respect to  $\mathcal{S}$ coincide. $\hfill\square$

The above corollary implies that we can take the factor group $\D_{\TT}/\mathcal{S}$.  Each element in $\D_{\TT}/\mathcal{S}$
is a conjugacy class containing a representative with the last column $(0,\dots,0,1)^T$.
Next, we note that the principal $5\times5$ sub-matrix of the matrix $D_2$ from \eqref{eq:T2-2} of the
transformation \eqref{eq:KP2019} is equal to that of the Kummer's transformation \eqref{eq:Kummer3F2}.
This transformation together with the permutation group $P_3\times{P_2}$ representing the obvious invariance with respect 
to separate permutations of the upper and lower parameters generate the entire group of Thomae transformations \cite{BLS1987}.
Next, comparing the principal $5\times5$ sub-matrices of the further generators $D_1$, $D_3$, $D_4$ with the matrices
of the Thomae transformations found, for instance in \cite[Appendix~1]{RaoDobermann}, we see that all of them occur among the
elements of the group of the Thomae transformations.  Hence, it remains to apply Theorem~3.2 from \cite{BLS1987} asserting 
that the group of the Thomae transformations is isomorphic the 120-element symmetric group $P_5$ of permutations on five symbols.  
Isomorphism is given by a linear change of variables seen in \eqref{eq:BLS}. Hence, our final result is the following theorem.

\begin{theorem}\label{th:Group-decomposed}
The group $\TT$ is isomorphic to $P_5\times\Z^5$.
\end{theorem}

As the entire group of the Thomae transformations for ${}_3F_{2}$ can be generated by the identity \eqref{eq:Kummer3F2}
and the permutation group $P_3\times{P_2}$, the above theorem implies that our entire group $\TT$
can be generated by the identity \eqref{eq:KP2019} (transformation $T_2$) and the top parameter shift transformation
$S_1$ together with the obvious symmetries $P_3\times{P_2}$. For example, the bottom parameter shift transformation can be obtained
as follows:
\begin{multline*}
(d-c,e-c,\psi,\psi+a,\psi+b)\stackrel{T_2^2}{\longmapsto}(a,b,c,d,e)\stackrel{S_1S_3^{-1}}{\longmapsto}(a+1,b,c-1,d,e)
\\
\stackrel{T_2^{-2}}{\longmapsto}(d-c+1,e-c+1,\psi,\psi+a+1,\psi+b)\stackrel{S_1S_2}{\longmapsto}(d-c,e-c,\psi,\psi+a+1,\psi+b).
\end{multline*}
Comparing the first and the last terms in this chain we see that we got  the bottom parameter shift transformation $S_4$ using only
$T_2$ and top shift transformations $S_1$, $S_2$, $S_3$ obtained from $S_1$ by permuting top parameters.

Theorem~\ref{th:Group-decomposed} further implies that there is a straightforward algorithm for computing any transformation 
from the group $\TT$. Details are given in the Appendix to this paper.
\end{section}

\begin{section}{Related ${}_3F_{2}$ transformation}

The proof of Theorem~\ref{th:Group} shows that each transformation $T\in\T$ is associated with the system \eqref{eq:sistem} of two ${}_3F_{2}$ transformations.  This system leads immediately to the following corollary.

\begin{proposition}\label{cr:Phi}
Each transformation $T\sim\{\varepsilon,M(\rr),\lambda(\rr),\alpha(\rr),\beta(\rr)\}\in\T$ induces a transformation for the ratio
$$
\Psi(\rr):=\frac{F_2(\rr)}{F_1(\rr)}=\frac{abc}{de}\frac{{}_{3}F_{2}\!\left(\begin{matrix}a+1,b+1,c+1\\d+1,e+1\end{matrix}\right)}
{{}_{3}F_{2}\!\left(\begin{matrix}a,b,c\\d,e\end{matrix}\right)}
=\frac{d}{dx}\log{}_{3}F_{2}\!\left(\left.\begin{matrix}a,b,c\\d,e\end{matrix}\right|x\right)_{|x=1}
$$
of the form
$$
\Psi(\rr)=\frac{\beta(\rr) \Psi(D\rr)+\lambda(\rr)}{\alpha(\rr)\Psi(D\rr)+\varepsilon}.
$$
\end{proposition}

Next, we observe that any two elements of $\T$ generate a three-term relation for ${}_3F_{2}$.
\begin{proposition}\label{prop:2}
For any two transformations from the group $\T$:
$T_1\sim\{\varepsilon_1, M_1(\rr)$,  $\lambda_1(\rr)$, $\alpha_1(\rr)$, $\beta_1(\rr)$, $ D_{1}\}$ and
$T_2\sim\{\varepsilon_2, M_2(\rr), \lambda_2(\rr),$ $\alpha_2(\rr),$ $\beta_2(\rr),\,D_{2}\}$
satisfying the condition $\alpha_2\beta_1-\alpha_1\beta_2\ne0$, the following identities hold
\begin{equation}\label{eq:threeterm1}
F_1(\rr)= M_1\frac{\alpha_2\beta_1\varepsilon_1-\alpha_1\alpha_2\lambda_1}{\alpha_2\beta_1-\alpha_1\beta_2} F_1(D_1\rr)+M_2\frac{\alpha_1\alpha_2\lambda_2-\alpha_1\beta_2\varepsilon_2}{\alpha_2\beta_1-\alpha_1\beta_2} F_1(D_2\rr)
\end{equation}
\emph{(}the dependence on $\rr$ is omitted for brevity\emph{)} and
\begin{equation}\label{eq:threeterm2}
F_2(\rr)=M_1\frac{\beta_1\beta_2\varepsilon_1- \alpha_1\beta_2 \lambda_1}{\alpha_2\beta_1 - \alpha_1\beta_2}F_1(D_1\rr) +M_2\frac{\alpha_2 \beta_1 \lambda_2-\beta_1\beta_2 \varepsilon_2}{\alpha_2 \beta_1 - \alpha_1 \beta_2}F_1(D_2\rr),
\end{equation}
where, as before, $F_1(\rr)={}_{3}F_{2}\!\left(\begin{matrix}a,b,c\\d,e\end{matrix}\right)$,  $F_2(\rr)=(abc)/(de){}_{3}F_{2}\!\left(\begin{matrix}a+1,b+1,c+1\\d+1,e+1\end{matrix}\right)$.
\end{proposition}

\textbf{Proof.} Solving \eqref{eq:sistem} for each transformation we, in particular,
get the system of equations:
$$
\left\{\begin{array}{l}F_1(D_1\rr)=(\beta_1F_1(\rr) - \alpha_1 F_2(\rr))/( M_1(\beta_1\varepsilon_1-\alpha_1\lambda_1)),
\\
F_1(D_2\rr)=(\beta_2F_1(\rr) - \alpha_2 F_2(\rr))/( M_2(\beta_2\varepsilon_2-\alpha_2\lambda_2)).\end{array}\right.
$$
Solving the above system for $F_1(\rr)$, $F_2(\rr)$ we arrive at \eqref{eq:threeterm1}, \eqref{eq:threeterm2}.
$\hfill\square$

If the matrices $D_{1}$, $D_{2}$ contain no shifts (i.e. the last column is $(0,0,0,0,0,1)^{T}$), then they correspond 
to Thomae's relations, so that $F_1(D_{1}\rr)$, $F_1(D_{2}\rr)$ are equal to each other up to a factor of gamma type.  
In this case, identities \eqref{eq:threeterm1}, \eqref{eq:threeterm2} become two-term transformations. 
However, for non-zero shifts Proposition~\ref{prop:2} generates genuine three-term relations
for ${}_{3}F_{2}\!\left(a,b,c;d,e\right)$. For example, we obtain
\begin{multline}\label{eq:example1}
    {}_{3}F_{2}\!\left(\begin{matrix}a,b,c\\d,e\end{matrix}\right)
    \\
    =\frac{\Gamma(d+1)\Gamma(e)\Gamma(d+e-a-b-c)}{\Gamma(a+1)\Gamma(d+e-a-b)\Gamma(d+e-a-c)}
    {}_{3}F_{2}\!\left(\begin{matrix}d+e-a-b-c-1, d-a, e-a
    \\
    d+e-a-c, e+d-a-b\end{matrix}\right)
    \\
    +\frac{(a-d)(d-b)(d-c)}{d(1+d)e}{}_{3}F_{2}\!\left(\begin{matrix}a+1,b+1,c+1\\d+2,e+1\end{matrix}\right).
\end{multline}
An important subclass of these transformations are pure shifts (the principal $5\times5$ submatrices of $D_{1}$, $D_{2}$ 
are identity matrices).  This subclass comprises the so-called contiguous relations, studied recently in detail in 
\cite{EbisuIwasaki}. In particular, Theorem~1.1 from \cite{EbisuIwasaki} claims the existence of the unique rational 
functions $u(\rr)$, $v(\rr)$ such that
\begin{equation}\label{eq:ThreeTermShift}
     {}_{3}F_{2}\!\left(\begin{matrix}a,b,c\\d,e\end{matrix}\right)=u(\rr) {}_{3}F_{2}\!\left(\begin{matrix}a+k_1,b+k_2,c+k_1\\d+k_4,e+k_5\end{matrix}\right)+  v(\rr) {}_{3}F_{2}\!\left(\begin{matrix}a+m_1,b+m_2,c+m_3\\d+m_4,e+m_5\end{matrix}\right)
\end{equation}
for any two distinct non-zero integer vectors $(k_1,\,k_2,\,k_3,\,k_4,\,k_5)$,  $(m_1,\,m_2,\,m_3,\,m_4,\,m_5)$.  Furthermore, Ebisu and Iwasaki presented a rather explicit algorithm in \cite{EbisuIwasaki} for computing the functions $u(\rr)$, $v(\rr)$ for given shifts.  Proposition~\ref{prop:2} furnishes an alternative method for computing these functions. For its realization we provide a collection of  \textit{Mathematica} routines in the Appendix to this paper.
Our algorithm works as follows: first step is to calculate transformations $T_1, T_2\in\TT$ associated with the matrices
$$
D_1=\begin{bmatrix}
1&0&0&0&0&k_1\\
0&1&0&0&0&k_2\\
0&0&1&0&0&k_3\\
0&0&0&1&0&k_4\\
0&0&0&0&1&k_5\\
0&0&0&0&0&1
\end{bmatrix},
D_2=\begin{bmatrix}
1&0&0&0&0&m_1\\
0&1&0&0&0&m_2\\
0&0&1&0&0&m_3\\
0&0&0&1&0&m_4\\
0&0&0&0&1&m_5\\
0&0&0&0&0&1
\end{bmatrix}.
$$
To this end we simply iterate transformations  $S^{\pm}$, $S_{\pm}$ realizing the shifts by  $\pm1$ of the first and 
forth parameters, respectively, combining them with the necessary permutations of the upper and lower parameters.  
To calculate the resulting $\lambda$, $\alpha$ and  $\beta$ the composition rule from Theorem~\ref{th:Group} is used 
with the help of \textit{Mathematica} routine. Then it remains to apply formula \eqref{eq:threeterm1}. For example, we get:
\begin{multline}\label{eq:ExampleThreeTerm}
{}_{3}F_{2}\!\left(\begin{matrix}a,b,c\\d,e\end{matrix}\right)=\frac{d+e-a-b-c-1}{e} {}_{3}F_{2}\!\left(\begin{matrix}a+1,b+1,c+1\\d+1,e+1\end{matrix}\right)
\\
+\frac{(a-d)(d-b)(d-c)}{d(d+1)e} {}_{3}F_{2}\!\left(\begin{matrix}a+1,b+1,c+1\\d+2,e+1\end{matrix}\right).
\end{multline}
Note that identity \eqref{eq:example1} is obtained from \eqref{eq:ExampleThreeTerm} by an application of a Thomae 
relation to the first term on the right hand side.  In a similar fashion, contiguous relations and Thomae transformations 
generate all three-term relations from Proposition~\ref{prop:2}, induced by the elements of the the group $\TT$.  
We note that the relations covered by Proposition~\ref{prop:2} are different from the three-term relations 
for ${}_{3}F_{2}$ summarized by Bailey in \cite[Section~3.7]{Bailey} and studied from group-theoretic viewpoint  in 
\cite[Section~IV]{BLS1987}.  This can be seen for example by comparing the matrices \cite[(2.6c)]{BLS1987} with the matrices $D$
associated with $\TT$.

The system  \eqref{eq:sistem} follows from the representation \eqref{eq:partition} of  ${}_{4}F_{3}$ with one unit shift 
as a linear combination of two ${}_{3}F_{2}$ functions. However, formula \eqref{eq:partition} is just one example 
of such decomposition. The two propositions that below give many more ways to expand the ${}_{4}F_{3}$ with unit shift 
into linear combination of ${}_3F_{2}$.  Proposition~\ref{pr:break} is proved directly in terms of hypergeometric series manipulations 
as its results will be used below in Section~6 to prove Lemma~\ref{lm:generators} used to generate the group $\TT$.
\begin{proposition}\label{pr:break}
The following identities hold true:
\begin{equation}\label{eq:shifttop1}
{}_{3}F_{2}\left.\!\!\left(\begin{matrix}\alpha,b,c\\d,e\end{matrix}\right.\right)
+\gamma{}_{3}F_{2}\left.\!\!\left(\begin{matrix}\alpha-1,b,c\\d,e\end{matrix}\right.\right)
=(\gamma+1){}_{4}F_{3}\left.\!\!\left(\begin{matrix}\alpha-1,b,c,\xi+1\\d,e,\xi\end{matrix}\right.\right),
\end{equation}
where $\xi=(\gamma+1)(\alpha-1)$;
\begin{equation}\label{eq:onetoponebottom1}
{}_{3}F_{2}\left.\!\!\left(\begin{matrix}\alpha,b,c\\d,e\end{matrix}\right.\right)
+\gamma{}_{3}F_{2}\left.\!\!\left(\begin{matrix}\alpha+1,b,c\\d+1,e\end{matrix}\right.\right)
=(\gamma+1){}_{4}F_{3}\left.\!\!\left(\begin{matrix}\alpha,b,c,\nu+1\\d+1,e,\nu\end{matrix}\right.\right),
\end{equation}
where $\nu={(\gamma+1)\alpha d}/(\gamma d+\alpha)$;  and
\begin{equation}\label{eq:allshiftedbutone1}
{}_{3}F_{2}\left.\!\!\left(\begin{matrix}\alpha,b,c\\d,e\end{matrix}\right.\right)+\gamma{}_{3}F_{2}\left.\!\!\left(\begin{matrix}\alpha,b+1,c+1\\d+1,e+1\end{matrix}\right.\right)
={}_{4}F_{3}\left.\!\!\left(\begin{matrix}\alpha-1,b,c,\lambda+1\\d,e,\lambda\end{matrix}\right.\right),
\end{equation}
where
$
\lambda=(\alpha-1)bc/(bc+\gamma d e).
$
\end{proposition}

\textbf{Proof.} We have
\begin{multline*}
{}_{3}F_{2}\left.\!\!\left(\begin{matrix}\alpha,b,c\\d,e\end{matrix}\right.\right)
+\gamma{}_{3}F_{2}\left.\!\!\left(\begin{matrix}\alpha-1,b,c\\d,e\end{matrix}\right.\right)
=
1+\gamma+\sum\limits_{n=1}^{\infty}\frac{(\alpha)_n(b)_n(c)_n+\gamma(\alpha-1)_n(b)_n(c)_n}{(d)_n(e)_nn!}
\\
=(1+\gamma)\left(1+\sum\limits_{n=1}^{\infty}\frac{(\alpha-1)_n(b)_n(c)_n}{(d)_n(e)_nn!}\left(1+\frac{n}{(\alpha-1)(\gamma+1)}\right)\right)\\=(\gamma+1){}_{4}F_{3}\left.\!\!\left(\begin{matrix}\alpha-1,b,c,\xi+1\\d,e,\xi\end{matrix}\right.\right),
\end{multline*}
where
$
\xi=(\gamma+1)(\alpha-1) $
 and we used $(\alpha)_n=(\alpha-1)_n(1+n/(\alpha-1))$.
Next,
\begin{multline*}
{}_{3}F_{2}\left.\!\!\left(\begin{matrix}\alpha,b,c\\d,e\end{matrix}\right.\right)
+\gamma{}_{3}F_{2}\left.\!\!\left(\begin{matrix}\alpha+1,b,c\\d+1,e\end{matrix}\right.\right)
=1+\gamma+\sum\limits_{n=1}^{\infty}\frac{(\alpha)_n(b)_n(c)_n}{(d+1)_n(e)_n n!}\left(1+\frac{ n}{d}+\gamma+\frac{\gamma n}{\alpha }\right)
\\
=(\gamma+1){}_{4}F_{3}\left.\!\!\left(\begin{matrix}\alpha,b,c,\nu+1\\d+1,e,\nu\end{matrix}\right.\right),
\end{multline*}
where
$\nu=(\gamma+1)\alpha d/(\gamma d+\alpha)$ and we used $(\alpha+1)_n=(\alpha)_n(1+n/\alpha)$.

Finally, using the obvious identities $(b)_n=b(b+1)_{n-1}$ and  $(\alpha)_n=(\alpha-1)_{n+1}/(\alpha-1)$  we get
\begin{multline*}
{}_{3}F_{2}\left.\!\!\left(\begin{matrix}\alpha,b,c\\d,e\end{matrix}\right.\right)
+\gamma{}_{3}F_{2}\left.\!\!\left(\begin{matrix}\alpha,b+1,c+1\\d+1,e+1\end{matrix}\right.\right)
=1+\sum\limits_{n=1}^{\infty}\frac{bc(\alpha)_{n-1}(b+1)_{n-1}(c+1)_{n-1}(\alpha+n-1)}{de(d+1)_{n-1}(e+1)_{n-1}n!}
\\
+\gamma{}_{3}F_{2}\left.\!\!\left(\begin{matrix}\alpha,b+1,c+1\\d+1,e+1\end{matrix}\right.\right)= 1+
 \sum\limits_{n=0}^{\infty}\frac{(\alpha)_{n}(b+1)_{n}(c+1)_{n}}{(d+1)_{n}(e+1)_{n}n!}\left(\frac{bc(\alpha+n)}{de(n+1)}+\gamma\right)
 \\
 =1+\sum\limits_{n=0}^{\infty}\frac{(\alpha-1)_{n+1}(b)_{n+1}(c)_{n+1}}{(d)_{n+1}(e)_{n+1}(n+1)!}\frac{de}{bc(\alpha-1)}\left(\frac{bc(\alpha+n)}{de}+\gamma(n+1)\right)=\\
 1+
 \sum\limits_{n=1}^{\infty}\frac{(\alpha-1)_{n}(b)_{n}(c)_{n}}{(d)_{n}(e)_{n}n!}\left(1+n\frac{bc+\gamma de }{(\alpha-1)bc}\right)={}_{4}F_{3}\left.\!\!\left(\begin{matrix}\alpha-1,b,c,\lambda+1\\d,e,\lambda\end{matrix}\right.\right),
 \end{multline*}
where $\lambda={(\alpha-1)bc}/{(bc+\gamma{de})}$.$\hfill\square$

Other ways to represent ${}_{4}F_{3}$ with one unit shift as a linear combination of ${}_{3}F_{2}$ are found 
by substituting  \eqref{eq:threeterm1} and \eqref{eq:threeterm2}  into  \eqref{eq:partition}.  
This is done in the following  proposition.
\begin{proposition}\label{prop:3}
Any two transformations from the group $\T$:
$T_1\sim\{\varepsilon_1, M_1(\rr)$, $\lambda_1(\rr),$  $\alpha_1(\rr),$ $\beta_1(\rr),\,$ $D_1\}$ and $T_2\sim\{\varepsilon_2,$ $M_2(\rr)$, $\lambda_2(\rr),$ $\alpha_2(\rr),$ $\beta_2(\rr),$ $ {D}_2\}$ satisfying the  condition $\alpha_2\beta_1-\alpha_1\beta_2\neq 0$ \emph{(}for brevity we omit the dependence on $\rr$ in the parameters\emph{)} induce the decomposition
\begin{equation}\label{eq:threeterm4}
{}_{4}F_3\!\!\left(\begin{matrix}a,b,c,f+1\\d,e,f\end{matrix}\right)= M_1\frac{\beta_1\varepsilon_1-\alpha_1\lambda_1}{\alpha_2\beta_1-\alpha_1\beta_2}\left(\alpha_2+\frac{\beta_2}{f}\right) F_1(D_1\rr)+M_2\frac{\alpha_2\lambda_2-\beta_2\varepsilon_2}{\alpha_2\beta_1-\alpha_1\beta_2}\left(\alpha_1+\frac{\beta_1}{f}\right) F_1(D_2\rr),
\end{equation}
where $F_1(\rr)={}_{3}F_{2}\!\left(\begin{matrix}a,b,c\\d,e\end{matrix}\right)$.
\end{proposition}
Let us exemplify \eqref{eq:threeterm4} with the following two decompositions:
\begin{multline*}
{}_{4}F_3\!\!\left(\begin{matrix}a,b,c,f+1\\d,e,f\end{matrix}\right)=\left(\frac{d+e-a-b-c-1}{e}+\frac{abc}{def}\right) {}_{3}F_2\!\!\left(\begin{matrix}a+1,b+1,c+1\\d+1,e+1\end{matrix}\right)
\\
+\frac{(a-d)(d-b)(d-c)}{ed(1+d)}{}_{3}F_2\!\!\left(\begin{matrix}a+1,b+1,c+1\\d+2,e+1\end{matrix}\right)
\end{multline*}
and
\begin{equation*}
{}_{4}F_3\!\!\left(\begin{matrix}a,b,c,f+1\\d,e,f\end{matrix}\right)=A{}_{3}F_2\!\!\left(\begin{matrix}a+1,b,c\\d,e\end{matrix}\right)+B{}_{3}F_2\!\!\left(\begin{matrix}a+1,b+1,c+1\\d+2,e+1\end{matrix}\right),
\end{equation*}
where
$$
A=1+\frac{bc(f-a)}{f(b(d-c)-d(d+e-a-c-1))},
~~
B=\frac{bc(a-d)(b-d)(c-d)(f-a)}{def(1+d)(b(c-d)+d(d+e-a-c-1))}.
$$
\end{section}

\begin{section}{Summation formulas}
In \cite[(45)]{KPDegenerate} we established the following summation formula
\begin{subequations}
\begin{equation}\label{eq:new4F3sum}
{}_{4}F_{3}\!\left(\begin{matrix}a,b, c, f+1\\d, e, f\end{matrix}\right)
=\frac{\Gamma(d)\Gamma(e)}{\Gamma(a+1)\Gamma(b+1)\Gamma(c+1)},
\end{equation}
valid if
\begin{equation}\label{eq:new4F3sum-restrictions}
e_1(d,e)-e_1(a,b,c)=2~~\text{and}~~f=\frac{e_3(a,b,c)}{e_2(a,b,c)-e_2(1-d,1-e)},
\end{equation}
\end{subequations}
where $e_k(\cdot)$ denotes the $k$-th elementary symmetric polynomial.   Now, if we apply any transformation of the form
\eqref{eq:generalform} and impose the above restrictions on the parameters on the right hand side, we obtain
\begin{subequations}
\begin{equation}\label{eq:4F3sum}
{}_{4}F_{3}\!\left(\begin{matrix}a,b, c, f+1\\d, e, f\end{matrix}\right)
=M(\rr)\frac{\varepsilon{f}+\lambda(\rr)}{f}F(\qq,\eta)=\frac{M(\rr)(\varepsilon{f}+\lambda(\rr))\Gamma(q_4)\Gamma(q_5)}{f\Gamma(q_1+1)\Gamma(q_2+1)\Gamma(q_3+1)},
\end{equation}
where $(q_1,q_2,q_3,q_4,q_5,1)=D\rr$, and the conditions $e_1(q_4,q_5)-e_1(q_1,q_2,q_3)=2$ and
$$
\eta=\frac{\varepsilon{f}+\lambda(\rr)}{\alpha(\rr)f+\beta(\rr)}=\frac{e_3(q_1,q_2,q_3)}{e_2(q_1,q_2,q_3)-e_2(1-q_4,1-q_5)}
$$
must hold.  Expressing $f$ these are equivalent to 
\begin{equation}\label{eq:4F3sum-restrictions}
e_1(q_4,q_5)-e_1(q_1,q_2,q_3)=2~\text{and}~
f=\frac{\lambda(\rr)(e_2(q_1,q_2,q_3)-e_2(1-q_4,1-q_5))-\beta(\rr)e_3(q_1,q_2,q_3)}
{\alpha(\rr)e_3(q_1,q_2,q_3)-\varepsilon(e_2(q_1,q_2,q_3)-e_2(1-q_4,1-q_5))}.
\end{equation}
\end{subequations}
As $q_{i}=q_i(a,b,c,d,e)$, $i=1,\ldots,5$, are linear functions we arrive at the following proposition:
\begin{proposition}\label{prop:summation}
Each transformation $T\in\T$ as characterized by the collection $\{\varepsilon, M(\rr),$ $\lambda(\rr),$ $ \alpha(\rr), $ $\beta(\rr),$ $ D\}$
corresponds to a summation formula \eqref{eq:4F3sum} valid under restrictions \eqref{eq:4F3sum-restrictions} with
$(q_1,\ldots,q_5,1)=D\rr$.
\end{proposition}
We will illustrate Proposition~\ref{prop:summation} by applying it to transformation \eqref{eq:T2}.
First condition in \eqref{eq:4F3sum-restrictions} becomes $e=c+2$. In view of this condition
formula \eqref{eq:4F3sum} takes the form
$$
{}_{4}F_{3}\!\left(\begin{matrix}a,b, c, f+1\\d, c+2, f\end{matrix}\right)
=\frac{(c + 1)\Gamma(d)
 \Gamma(d - a - b + 
   2)(f\psi + c(a + b - d))}{\Gamma(d - a + 1)\Gamma(d - b + 1)f\psi},
$$
where $\psi=d-a-b+1$ and, by the second condition in \eqref{eq:4F3sum-restrictions},
$$
f=-\frac{c(a + b - d)}{\psi} + \frac{(d - a)(d - b) c}{\psi((d - a)(d - b + c) + (d - b)
        c + (d - 1)(a + b - d - c - 1))}.
$$
Further examples will be given in \cite{CK2020}.
\end{section}

\begin{section}{Proof of Lemma~\ref{lm:generators}}
Write identity \eqref{eq:partition} in expanded form
\begin{equation}\label{4F3-decomposition}
{}_{4}F_{3}\!\!\left(\begin{matrix}a,b,c,f+1\\d,e,f\end{matrix}\right)
={}_{3}F_{2}\!\!\left(\begin{matrix}a,b,c\\d,e\end{matrix}\right)+\frac{abc}{fde}{}_{3}F_{2}\!\!\left(\begin{matrix}a+1,b+1,c+1\\d+1,e+1\end{matrix}\right).
\end{equation}
Applying Thomae's transformation \eqref{eq:Thomae3F2} to both ${}_{3}F_{2}$ functions on the right hand side, 
we get ($\psi=d+e-a-b-c-1$):
\begin{multline*}
{}_{4}F_{3}\!\!\left(\begin{matrix}a,b,c,f+1\\d,e,f\end{matrix}\right)
=\frac{\Gamma(\psi+1)\Gamma(d)\Gamma(e)}{\Gamma(a)\Gamma(\psi+b+1)\Gamma(\psi+c+1)}\times
\\
\left[{}_{3}F_{2}\!\!\left(\begin{matrix}\psi+1,d-a,e-a\\\psi+b+1,\psi+c+1\end{matrix}\right)
+\frac{bc}{f\psi}{}_{3}F_{2}\!\!\left(\begin{matrix}\psi,d-a,e-a\\\psi+b+1,\psi+c+1\end{matrix}\right)\right].
\end{multline*}

Now we employ Proposition~\ref{pr:break}.  Application of formula \eqref{eq:shifttop1} to the linear combination in brackets yields
$$
{}_{4}F_{3}\!\!\left(\begin{matrix}a,b,c,f+1\\d,e,f\end{matrix}\right)
=\frac{(f\psi+bc)\Gamma(\psi)\Gamma(d)\Gamma(e)}{f\Gamma(a)\Gamma(\psi+b+1)\Gamma(\psi+c+1)}
{}_{4}F_{3}\!\!\left(\begin{matrix}\psi,d-a,e-a,\eta+1\\\psi+b+1,\psi+c+1,\eta\end{matrix}\right),
$$
where $\eta=\psi+bc/f$.   This proves transformation given by \eqref{eq:generator1}.

In a similar fashion, if we apply the Kummer transformation \eqref{eq:Kummer3F2} to ${}_{3}F_{2}$ on the right  hand side
of \eqref{4F3-decomposition} we get:
\begin{multline*}
{}_{4}F_{3}\!\!\left(\begin{matrix}a,b,c,f+1\\d,e,f\end{matrix}\right)
=\frac{\Gamma(\psi+1)\Gamma(d)}{\Gamma(d-a)\Gamma(\psi+a+1)}\times
\\
\left[{}_{3}F_{2}\!\!\left(\begin{matrix}a,e-b,e-c\\e,\psi+a+1\end{matrix}\right)
+\frac{abc}{fe\psi}{}_{3}F_{2}\!\!\left(\begin{matrix}a+1,e-b,e-c\\e+1,\psi+a+1\end{matrix}\right)\right].
\end{multline*}

Applying the relation \eqref{eq:onetoponebottom1} to the linear combination in brackets we then obtain
$$
{}_{4}F_{3}\!\!\left(\begin{matrix}a,b,c,f+1\\d,e,f\end{matrix}\right)
=\frac{(abc+fe\psi)\Gamma(\psi)\Gamma(d)}{fe\Gamma(d-a)\Gamma(\psi+a+1)}
{}_{4}F_{3}\!\!\left(\begin{matrix}a,e-b,e-c,\lambda+1\\e+1,\psi+a+1,\lambda\end{matrix}\right),
$$
where
$$
\lambda=\frac{abc+fe\psi}{bc+f\psi}.
$$
This proves transformation \eqref{eq:generator3}.
\end{section}

\begin{section}{Appendix}
In this appendix we will display the  explicit form of the main building blocks needed for calculating the elements of 
the group $\TT$. Just like with Thomae's transformations \cite[Appendix~1]{RaoDobermann} we have 10 different identities 
with zero shifts.  They are obtained as follows: permuting $a\leftrightarrow{b}$ and $a\leftrightarrow{c}$ in 
formula \eqref{eq:generator1} we get three transformations while $a\leftrightarrow{b}$, $a\leftrightarrow{c}$ 
and $d\leftrightarrow{e}$ in \eqref{eq:KP2019} leads to six more transformations. Adding the identity transformation 
we arrive at 10 ''Thomae-like'' zero-shift  transformations for ${}_4F_3$ with unit shift.  The entire 120 element 
subgroup of   ''Thomae-like'' zero-shift  transformations is obtained by the obvious 12 permutations of three top and 
two bottom parameters on the right hand side of each of the 10 transformations described above.

All further transformations are obtained by consecutive application of the four shifting 
trans\-for\-ma\-ti\-ons $S^{\pm}$, $S_{\pm}$ and permutations of top and bottom parameters to the 120 
transformations described above.  Transformation $S^{+}$ shifting the top parameter $a$ by $+1$ (denoted by $S_{1}$ in Section~3) 
is given by \eqref{eq:killuppershift}.  Combining parameters it can be written as:
\begin{equation}\label{eq:killuppershift1}
{}_{4}F_3\!\!\left(\begin{matrix}a,b,c,f+1\\d,e,f\end{matrix}\right)
=\left(1-\frac{bc}{(d-a-1)(e-a-1)}\right)\left(1+\frac{\lambda}{f}\right)
{}_{4}F_3\!\!\left(\begin{matrix}a+1,b,c,\eta+1\\d,e,\eta\end{matrix}\right),
\end{equation}
where
$$
\lambda=\frac{abc}{a(2+a-d-e)-bc+(d-1)(e-1)},~~\eta=\frac{abc+((a+1)(a+1-d-e)-bc+de)f}
{(a-f)(2+a+b+c-d-e)}.
$$
Its inverse $S^{-}$ is given by:
\begin{equation}\label{eq:killuppershift-1}
{}_{4}F_3\!\!\left(\begin{matrix}a,b,c,f+1\\d,e,f\end{matrix}\right)
=\left(1+\frac{bc}{\psi{f}}\right)
{}_{4}F_3\!\!\left(\begin{matrix}a-1,b,c,\eta+1\\d,e,\eta\end{matrix}\right),
\end{equation}
where
$$
\eta=\frac{(a-1)(bc+\psi{f})}{a(d+e-a)+bc-de+\psi{f}}.
$$
The transformation $S_{+}$ shifting the bottom parameter $d$ by $+1$ (denoted by $S_{4}$ in Section~3) 
is given by \eqref{eq:killbottom}.  It can be written more compactly as
\begin{equation}\label{eq:killbottom1}
{}_{4}F_3\!\!\left(\begin{matrix}a,b,c,f+1\\d,e,f\end{matrix}\right)
=\frac{abc+\psi{d}f}{\psi{d}f}
{}_{4}F_3\!\!\left(\begin{matrix}a,b,c,\eta+1\\d+1,e,\eta\end{matrix}\right),
\end{equation}
where $\psi=e+d-a-b-c-1$ and
$$
\eta=\frac{abc+\psi{d}f}{d(d-a-b-c)+ab+ac+bc+\psi{f}}.
$$
Finally, its inverse transformation $S_{-}$ shifting a bottom parameter by $-1$ has the form
\begin{equation}\label{eq:killbottom-1}
{}_{4}F_3\!\!\left(\begin{matrix}a,b,c,f+1\\d,e,f\end{matrix}\right)
=\frac{[((d-b-1)(d-c-1)-a(d-b-c-1))f-abc](d-1)}{(d-a-1)(d-b-1)(d-c-1)f}
{}_{4}F_3\!\!\left(\begin{matrix}a,b,c,\eta+1\\d-1,e,\eta\end{matrix}\right),
\end{equation}
where
$$
\eta=\frac{abc+[(1-d)(d-a-b-c-1)-ab-ac-bc]f}{(d+e-a-b-c-2)(f-d+1)}.
$$

In the remaining part of the Appendix we present several  \textit{Wolfram Mathematica}\textsuperscript{\tiny\textregistered} 
routines intended for dealing with the group $\T$ together with an example of their use.  Listing~\ref{lst:comp} contains the 
function $\mathrm{CMPS}[T_1,T_2]$ that takes as input two transformations $T_1$, $T_2$ and computes their 
composition $T_2\circ{T_1}$. The form in which the parameters $\varepsilon_i$, $M_i$, $\lambda_i$, $\alpha_i$, $\beta_i$ and $D_i$, 
$i=1,2$, should be supplied can be seen from the example in Listing~\ref{lst:example}.  Similarly, Listing~\ref{lst:inv} 
contains the function $\mathrm{INV}[T]$ that  computes the inverse of a given transformation $T$.  The output provided 
by $\mathrm{CMPS}$ and $\mathrm{INV}$ can be printed in a easily readable form using the function $\mathrm{PRN}[T]$ given in 
Listing~\ref{lst:prn}.  The same Listing~\ref{lst:prn} contains the function $\mathrm{INPT}[T]$ that converts the output form 
of the functions $\mathrm{CMPS}$ and $\mathrm{INV}$ into the input form of the same functions, so that further compositions 
or inverses could be computed from such output. For numerical verification of the outputs of $\mathrm{CMPS}$ and $\mathrm{INV}$ 
the function $\mathrm{RHS}[T]$ presented in Listing~\ref{lst:rhs} converts these outputs into the expression that can 
be evaluated by the  \textit{Mathematica} function $\mathrm{N}[...]$ after the parameters have been assigned some numerical values, 
see an example at the end of Listing~\ref{lst:example}.

\begin{lstlisting}[language=Mathematica,caption={Composition},label={lst:comp}]
   CMPS[T1_, T2_]:=Module[{eps1=T1[[1]], M1=T1[[2]], lam1=T1[[3]], alpha1=T1[[4]],
   beta1=T1[[5]], D1=T1[[6]], eps2=T2[[1]], M2=T2[[2]], lam2=T2[[3]], alpha2=T2[[4]], beta2=T2[[5]],
   D2=T2[[6]], R={{a},{b},{c},{d},{e},{1}}}, RR=Flatten[Drop[R,{6}]];
   If[Simplify[eps1*eps2+alpha1@@RR*lam2@@Flatten[Drop[D1.R,{6}]]]===0,
   {0, FullSimplify[M1@@RR*M2@@Flatten[Drop[D1.R,{6}]]*(eps2*lam1@@RR+lam2@@Flatten[Drop[D1.R, {6}]]*beta1@@RR)], 1,
   Simplify[(eps1*alpha2 @@ Flatten[Drop[D1.R,{6}]]+alpha1@@RR*beta2@@Flatten[Drop[D1.R, {6}]])/
   (eps2*lam1@@RR+lam2@@Flatten[Drop[D1.R,{6}]]*beta1@@RR)], Simplify[(lam1@@RR*alpha2@@Flatten[Drop[D1.R,{6}]]
   +beta1@@RR*beta2@@Flatten[Drop[D1.R, {6}]])/(eps2*lam1@@RR+lam2@@Flatten[Drop[D1.R,{6}]]*beta1@@RR)], D2.D1},
   {1, FullSimplify[M1@@RR*M2@@Flatten[Drop[D1.R,{6}]]*(eps1*eps2+alpha1@@RR*lam2@@Flatten[Drop[D1.R, {6}]])],
   Simplify[(eps2*lam1@@RR+lam2@@Flatten[Drop[D1.R,{6}]]*beta1@@RR)/
   (eps1*eps2+alpha1@@RR*lam2@@Flatten[Drop[D1.R,{6}]])], Simplify[(eps1*alpha2@@Flatten[Drop[D1.R, {6}]]
   +alpha1@@RR*beta2@@Flatten[Drop[D1.R, {6}]])/(eps1*eps2+alpha1@@RR*lam2@@Flatten[Drop[D1.R, {6}]])],
   Simplify[(lam1@@RR*alpha2@@Flatten[Drop[D1.R,{6}]]+beta1@@RR*beta2 @@Flatten[Drop[D1.R, {6}]])/
   (eps1*eps2+alpha1@@RR*lam2@@Flatten[Drop[D1.R,{6}]])], D2.D1}]]
  \end{lstlisting}

\begin{lstlisting}[language=Mathematica,caption={Inversion},label={lst:inv}]
INV[TT_]:=Module[{eps=TT[[1]], M=TT[[2]], lam=TT[[3]], alpha=TT[[4]], beta=TT[[5]], D=TT[[6]],
   R={{a}, {b}, {c}, {d}, {e}, {1}}}, RR=Flatten[Drop[R, {6}]];
   If[Simplify[beta@@RR]===0, {0,FullSimplify[1/M@@Flatten[Drop[Inverse[D].R,{6}]]/
   alpha@@Flatten[Drop[Inverse[D].R, {6}]]], 1, Simplify[alpha@@Flatten[Drop[Inverse[D].R,{6}]]/
   lam@@Flatten[Drop[Inverse[D].R, {6}]]], -eps/lam@@Flatten[Drop[Inverse[D].R, {6}]], Inverse[D]},
   {1, FullSimplify[beta@@Flatten[Drop[Inverse[D].R, {6}]]/(M@@Flatten[Drop[Inverse[D].R,{6}]]*
   (eps*beta@@Flatten[Drop[Inverse[D].R, {6}]]-lam@@Flatten[Drop[Inverse[D].R, {6}]]*
   alpha@@Flatten[Drop[Inverse[D].R, {6}]]))], Simplify[-lam@@Flatten[Drop[Inverse[D].R,{6}]]/
   beta@@Flatten[Drop[Inverse[D].R, {6}]]], Simplify[-alpha@@Flatten[Drop[Inverse[D].R,{6}]]/
   beta@@Flatten[Drop[Inverse[D].R, {6}]]], Simplify[eps/beta @@ Flatten[Drop[Inverse[D].R,{6}]]], Inverse[D]}]]
\end{lstlisting}

\begin{lstlisting}[language=Mathematica,caption={Conversion into input form and printing},label={lst:prn}]
exprToFunction[expr_, vars_]:=ToExpression[ToString[FullForm[expr]/.MapIndexed[#1 -> Slot @@ #2 &, vars]]<>"&"];

INPT[TT_]:=List[TT[[1]], exprToFunction[TT[[2]], {a, b, c, d, e}],
  exprToFunction[TT[[3]], {a, b, c, d, e}], exprToFunction[TT[[4]], {a, b, c, d, e}],
  exprToFunction[TT[[5]], {a, b, c, d, e}], TT[[6]]]

ETA[TT_]:=Collect[Numerator[Together[(TT[[1]]*f+TT[[3]])/(TT[[4]]*f+TT[[5]])]], f]/
  Collect[Denominator[Together[(TT[[1]]*f + TT[[3]])/(TT[[4]]*f + TT[[5]])]], f]

PRN[TT_]:=Module[{}, Print["epsilon=", TT[[1]]];
   Print["M=", FullSimplify[TT[[2]]]];
   Print["Lambda=", FullSimplify[TT[[3]]]]; Print["alpha=", TT[[4]]];
   Print["Beta=", TT[[5]]];
   Print["Parameters=", Flatten[Drop[TT[[6]].{{a}, {b}, {c}, {d}, {e}, {1}}, {6}]]];
   Print["eta=", ETA[TT]]];
\end{lstlisting}

\begin{lstlisting}[language=Mathematica,caption={Conversion into computable form},label={lst:rhs}]
 RHS[TT_]:=Simplify[TT[[2]]*(TT[[1]]*f + TT[[3]])/f]*
 HypergeometricPFQ[Join[Flatten[Drop[TT[[6]].{{a},{b},{c},{d},{e},{1}},{6}]][[1;;3]], {ETA[TT]+1}],
  Join[Flatten[Drop[TT[[6]].{{a},{b},{c},{d},{e},{1}},{6}]][[4;;5]], {ETA[TT]}],1]
\end{lstlisting}

\begin{lstlisting}[language=Mathematica,caption={Example of use},label={lst:example}]
(*Definition of the first transformation*)
eps1=1; M1[a_,b_,c_,d_,e_]:=Gamma[d+e-a-b-c]*Gamma[d]*Gamma[e]/Gamma[a]/Gamma[d+e-a-c]/Gamma[d+e-a-b];
lam1[a_,b_,c_,d_,e_]:=b*c/(d+e-a-b-c-1); alpha1[a_,b_,c_,d_,e_]:=1/(d+e-a-b-c-1);
beta1[a_, b_, c_, d_,e_]:=0; D1={{-1,-1,-1,1,1,-1}, {-1,0,0,1,0,0}, {-1,0,0,0,1,0}, {-1,0,-1,1,1,0},
{-1,-1,0,1,1,0}, {0,0,0,0,0,1}};

(*Definition of the second transformation*)
eps2=1; M2[a_,b_,c_,d_,e_]:=Gamma[d+e-a-b-c]*Gamma[e]/Gamma[d+e-a-b]/Gamma[e-c];
lam2[a_,b_,c_,d_,e_]:=(a+b-d)*c/(d+e-a-b-c-1); alpha2[a_,b_,c_,d_,e_]:=0;
beta2[a_,b_,c_,d_,e_]:=(e-c-1)/(d+e-a-b-c-1); D2={{-1,0,0,1,0,0}, {0,-1,0,1,0,0}, {0,0,1,0,0,0},
{0,0,0,1,0,0}, {-1,-1,0,1,1,0}, {0,0,0,0,0,1}};

(*composition T1T2*)
T1T2=CMPS[{eps1, M1, lam1, alpha1, beta1, D1}, {eps2, M2, lam2, alpha2, beta2, D2}];

(*Inverse of T1*))
T1INV = INV[{eps1, M1, lam1, alpha1, beta1, D1}];

(*Printing the parameters of T1T2*)
PRN[T1T2]
epsilon=1
M=-(((a c+(1+b-d) e) Gamma[d] Gamma[-1-a-b-c+d+e])/(e Gamma[-b+d] Gamma[-a-c+d+e]))
Lambda=-((a b c)/(a c+(1+b-d)e))
alpha=(1+b-d)/(a c+e+b e-d e)
Beta=0
Parameters={1+b,-c+e,-a+e,-a-c+d+e,1+e}
eta=(-abc+(ac+e+be-de)f)/((1+b-d)f)

(*Taking composition if the results of previous operations*))
NEW=CMPS[INPT[T1T2], INPT[T1INV]];

(*Numerical verification o the transformation NEW using RHS[...]*)
a=1+2/3; b=-13/17+2; c=3/7; d=5/11; e=5+44/17; f=12/13;
In[51]:= N[HypergeometricPFQ[{a, b, c, f + 1}, {d, e, f}, 1], 15]
Out[51]= 2.22268615827388
In[52]:= N[RHS[NEW], 15]
Out[52]= 2.22268615827388
\end{lstlisting}

\end{section}

\paragraph{Acknowledgements.} The second author is supported by the Ministry of Science and Higher Education of the Russian Federation (supplementary agreement No.075-02-2020-1482-1 of April 21, 2020).

\end{document}